\documentclass[doublespacing]{elsart}
\usepackage{amssymb}
\usepackage[english,francais]{babel}
\usepackage[latin1]{inputenc}
\usepackage{amsmath}
\newtheorem{theorem}{Theorem}[section]
\newtheorem{lemma}[theorem]{Lemma}
\newtheorem{e-proposition}[theorem]{Proposition}
\newtheorem{corollary}[theorem]{Corollary}
\newtheorem{e-definition}[theorem]{Definition\rm}

\newtheorem{theoreme}{Th\'eor\`eme}[section]

\newtheorem{proposition}[theoreme]{Proposition}

\def\og{\leavevmode\raise.3ex\hbox{$\scriptscriptstyle\langle\!\langle$~}}
\def\fg{\leavevmode\raise.3ex\hbox{~$\!\scriptscriptstyle\,\rangle\!\rangle$}}
\begin{document}
\centerline{}
\begin{frontmatter}
\selectlanguage{english}
\title{On  a projectively invariant distance on  Einstein Finsler spaces}
\selectlanguage{english}
\author{M. SEPASI and B. BIDABAD\corauthref{cor1}\thanksref{label2}}
 \thanks[label2]{The second author is partially supported by INSF no,890000676.}
 \corauth[cor1]{Corresponding author.}
\ead{m\_sepasi@aut.ac.ir;bidabad@aut.ac.ir;bbidabad@iupui.edu}
\address{Department of Mathematics, Amirkabir University of Technology \\(Tehran Polytechnic), Tehran 15914, Iran.}
\begin{abstract}
In this work  an intrinsic projectively invariant distance is used to establish a new approach to the study of projective geometry in Finsler space. It is shown that the projectively invariant distance previously defined is a constant multiple of the Finsler distance in certain case.
 As a consequence,   two projectively related complete Einstein Finsler spaces with constant negative scalar curvature are homothetic. Evidently, this will be true for Finsler spaces of constant flag curvature as well.
 \vskip
0.5\baselineskip
\noindent{\bf R\'esum\'e} \vskip 0.5\baselineskip \noindent {\bf Sur une distance projectivement invariante en espace d'Einstein-Finsler.}
\selectlanguage{francais}
Dans ce travail, une distance intrins\`{e}que projectivement invariante est utilis\'{e} pour \'{e}tablir une nouvelle approche pour l'\'{e}tude de la g\'{e}om\'{e}trie projective dans l'espace de Finsler. Il est montr\'{e} que la distance projectivement invariante d\'{e}finie pr\'{e}c\'{e}demment est un multiple constant de la distance finslerienne dans certain cas. Par cons\'{e}quent, deux espaces d'Einstein-Finsler  compl\'{e}tes  \`{a}  courbure scalaire constante n\'{e}gative sont homoth\'{e}tiques. Évidemment, cela sera vrai aussi  pour les espaces de Finsler \`{a} courbure sectionelle constante.
 \end{abstract}
\end{frontmatter}
\selectlanguage{english}
{\small\emph{Mathematics Subject Classification}: Primary 53C60; Secondary 58B20.}
\section*{Introduction}
Two regular metrics on a manifold are said to be pointwise projectively related if they have the same geodesics as the point sets. Two regular metric spaces are said to
be projectively related if there is a diffeomorphism between them such that the pull-back of one metric is
pointwise projective to the other.  Let $\gamma$ be a geodesic of a metric space.
In general, the parameter $t$  of $\gamma(t)$ does not remain invariant under the projective changes. There is a unique parameter up to linear fractional  transformations which is projectively invariant. This parameter is referred to, in the literature, as \emph{projective parameter}. See \cite{Be,SB} for a survey. The projective parameter together with the Funk metric  is used to establish a projectively invariant pseudo-distance in Finsler spaces. Next a comparison theorem on Ricci curvatures shows that this pseudo-distance is a distance.
  The  Ricci tensor was introduced in Riemannian spaces 1904 by G. Ricci and nine years later  was used to formulate Einstein's  theory of gravitation \cite{Bo}. In the presentwork we use the notion of Ricci curvature introduced by Akbar-Zadeh, cf. \cite{A}.
   Hence a Finsler metric is said to be Einstein if the Ricci scalar Ric is a function of $x$ alone. Equivalently $Ric_{ij}=Ric(x)g_{ij}$.

 Without pretending to be exhaustive, we bring a few  results related to our approach on Einstein-Finsler space. If M is simply connected endowed with a complete metric connection such that the symmetric part of the Ricci curvature of the associated symmetric connection is of Einstein type, that is, $R_{(ij)}=cg_{ij}$, where $c$ is a positive constant and if M admits a projective group leaving invariant the trace of torsion, then M is homeomorphic to an sphere, cf., \cite{AC}.
In \cite{S} Z. Shen found out  that  two pointwise projectively equivalent  Einstein Finsler metrics $F$ and $\bar{F}$ on an n-dimensional compact manifold $M$ have same  sign Einstein constants. In addition, if two pointwise projectively related
Einstein metrics are complete with negative Einstein constants, then one of them is a multiple of the other. Later in a joint work, he proved that if two projectively related Riemannian metrics $g$ and $\bar{g}$ on a manifold $M$ have Ricci curvatures satisfying $\bar{Ric} \leq Ric$ and $g$ is complete, then the projective change is affine.\cite{CS}. Recently, G. Yang generalized this comparison on Ricci curvatures of Finsler spaces and get some interesting results about length of geodesics  and completeness of the space \cite{Y}.
Here, inspired by  the Kobayashi's work \cite{K},  the projectively invariant distance   in complete Einstein spaces is studied and it is proved that   the intrinsic distance is a constant multiple of Finsler distance. Consequently the topology generated by the intrinsic distance coincide with that of Finslerian distance and in a new approach, we find out the known fact  that two projectively related complete Einstein Finsler spaces with constant negative Ricci scalar are homothetic, cf. \cite{S}.
\section{Preliminaries}
Let $M$ be an $n$-dimensional  $C^{\infty}$  connected manifold. Denote by $T_xM$ the tangent space at $x\in M$, and by $TM:={\cup}_{x\in M} T_xM$ the  bundle of tangent spaces. Each element of $TM$ has the form $(x,y)$, where $x\in M$ and $y\in T_xM$. The natural projection $\pi:TM\rightarrow M $, is given by $\pi (x,y):= x$.  The pull-back tangent bundle ${\pi}^{*} TM$ is a vector bundle over the slit tangent bundle $TM_0:=TM\backslash\{0\}$ for which the fiber ${\pi}^{*}_vTM$ at $v \in TM_0$ is just $T_xM$, where $\pi (v) = x$. \\
A (globally defined) Finsler structure on $M$ is a function $F: TM\rightarrow [0 , \infty) $ with the  properties; (I) Regularity: $F$ is $C^{\infty}$ on the entire slit tangent bundle $TM_0$; (II) Positive homogeneity: $F(x , \lambda y) = \lambda F(x , y)$ for all $\lambda > 0$; (III) Strong convexity: The $n \times n$ Hessian matrix $(g_{ij}):= ({[\frac{1}{2}F^2]}_{y^iy^j})$, is positive-definite at every point of $TM_0$.
For  any $y\in T_xM_0$, the Hessian $g_{ij}(y)$ induces an inner product $g_y$ in $T_xM$ by $g_y(u,v):=g_{ij}(y)u^iv^j$.
Let  $\gamma : [a , b] \rightarrow M$ be a piecewise $C^{\infty}$ curve on $(M,f)$ with the velocity $\frac{d\gamma}{dt} = \frac{d{\gamma}^i}{dt} \frac{\partial}{\partial x^i} \in T_{\gamma (t)} M$. The arc length parameter of $\gamma$ is given by $s(t) = \int_{t_0}^{t} F(\gamma , \frac{d\gamma}{dr}) dr$, and the integral length  is denoted by $L(\gamma):= \int_a^b F(\gamma , \frac{d\gamma}{dt}) dt $.
 For every $x_0$ , $x_1$ $\in M$, denote by $\Gamma (x_0 , x_1)$ the collection of all piecewise $C^{\infty}$ curves  $\gamma : [a , b] \rightarrow M$ with $\gamma (a) = x_0$ and $\gamma(b) = x_1$, and define a map $d_F : M \times M \rightarrow [0 , \infty)$ by $d_F(x_0 , x_1) := inf L(\alpha)$, where $\alpha \in \Gamma (x_0 , x_1)$. It can be shown that $d_F$ satisfies the first two axioms of a metric space. Namely, (I) $d_F(x_0 , x_1) \geq 0$ , where equality holds if and only if $x_0 = x_1$; (II) $d_F(x_0 , x_1) \leq d_F(x_0 , x_1) + d_F(x_1 , x_2)$.\\
We should remark that the distance function $d_F$ on a Finsler space does not have  the symmetry property.
If the Finsler structure F is absolutely homogeneous, that is $F(x,\lambda y)=\mid\lambda\mid F(x,y)$ for $\lambda \in \mathbb{R}$, then one also has the third axiom of a metric space,  (III) $d(x_0, x_1) = d(x_1,x_0)$. The manifold topology coincides with that
generated by the forward metric balls, $B_p^+(r) := \{x \in M : d_F(p , x) < r\}$. The latter assertion is also true for backward metric balls, $B_p^-(r)$, cf.,  \cite{BCS}.
Every Finsler metric $F$ induces a spray $\textbf{G}=y^i\frac{\partial}{\partial x^i}-2G^i(x,y)\frac{\partial}{\partial y^i}$ on $TM$, where $G^i(x,y):=\frac{1}{4} g^{il}\{[F^2]_{x^k y^l}y^k-[F^2]_{x^l}\}$.
$\textbf{G}$ is a globally defined vector field on $TM$.  Projection of a flow line of $\textbf{G}$
is called a geodesic on M. A curve $\gamma(t)$ on $M$ is a geodesic if and
only if in local coordinate it satisfies $\frac{d^2x^i}{ds^2}+2G^i(x(s),\frac{dx}{ds})=0$, where $s$ is the arc length parameter.
$F$ is said to be positively complete (resp. negatively complete), if any geodesic on
an open interval $(a, b)$ can be extended to a geodesic on $(a,\infty)$ (resp. $(-\infty, b)$).
$F$ is said to be complete if it is positively and negatively complete.
 For a vector $y \in T_xM_0$, the Riemann curvature
$\textbf{R}_y : T_xM \rightarrow T_xM$ is defined by $\textbf{R}_y(u)=R^i_k u^k \frac{\partial}{\partial x^i}$, where $R^i_k(y):=2\frac{\partial G^i}{\partial x^k}-\frac{{\partial}^2G^i}{\partial y^k \partial x^j}y^j+2G^j \frac{{\partial}^2G^i}{\partial y^k \partial y^j}-\frac{\partial G^i}{\partial y^j}\frac{\partial G^j}{\partial y^k}$.\\
For a two-dimensional plane $P \subset T_pM$ and a non-zero vector $y\in T_pM$, the flag curvature $\textbf{K}(P, y )$ is defined by $\textbf{K}(P,y):=\frac{g_y(u,\textbf{R}_y(u))}{g_y(y,y)g_y(u,u)-{g_y(y,u)}^2}$, where $P=$ span$\{y,u\}$. $F$ is said to be of \textit{scalar curvature} $\textbf{K}=\lambda(y)$ if for any $y\in T_pM$, the flag curvature $\textbf{K}(P,y)=\lambda(y)$ is independent of $P$ containing $y\in T_pM$. It is equivalent to the following system in a local coordinate system $(x^i,y^i)$ on                                                                               $TM$.
\begin{equation}\label{scalar cur}
R^i_k=\lambda F^2\{{\delta}^i_k-F^{-1}F_{y^k}y^i\}.
\end{equation}
If $\lambda$ is a constant, then $F$ is said to be of \textit{constant curvature}. The Ricci scalar  of $F$ is a positive zero homogeneous function in $y$ given by $Ric:=\frac{1}{F^2} R^i_i$. This is equivalent to say that
$Ric (x , y )$ depends on the direction of the flag pole $y$ but not its length. The Ricci tensor of a Finsler
metric $F$ is defined by $Ric_{ij}:=\{\frac{1}{2}R^k_k\}_{y^iy^j}$, cf., \cite{A}.
If $(M,F)$ is a Finsler space with constant flag curvature $\lambda$,  (\ref{scalar cur}) leads to
\begin{equation}\label{ric,k}
Ric=(n-1)\lambda, \qquad Ric_{ij}= (n-1) \lambda g_{ij}.
\end{equation}
A Finsler metric is said to be  an \textit{Einstein metric} if the Ricci scalar function is a function of $x$ alone, equivalently $Ric_{ij}=Ric(x)g_{ij}$.
\subsection{Projective parameter and Schwarzian derivative}
 A Finsler space $(M, F)$ is said to be \textit{projective }to another Finsler space $(M,\bar{F})$ as a set of points, if and only if there exists a one-positive homogeneous scalar field $P(x,y)$ on $TM$ satisfying ${\bar{G}}^i(x,y)= G^i(x,y)+P(x,y)y^i$. The scalar field $P(x,y)$ is called the \emph{projective factor} of the projective change under consideration.
In general, the parameter $t$ of a geodesic, does not remain invariant under projective change of metrics. It is well known, there  is a unique parameter up to the linear fractional  transformations which is projectively invariant. This parameter is referred to, in the literature, as \emph{ projective parameter}. The projective parameter, for a geodesic  $\gamma$ is given by
$\{\pi,s\}= \frac{2}{n-1} Ric_{jk}\frac{d{x}^j}{ds}\frac{d{x}^k}{ds}$, where the operator $\{.,.\}$ is the \emph{Schwarzian derivative}  defined for a $C^{\infty}$ real function $f$ on $\mathbb{R}$, and for $t\in \mathbb{R}$ by  $\{f,t\}=\frac{{f}^{'''}}{{f}^{'}}-3/2{(\frac{{f}^{''}}{{f}^{'}})}^2$, where $f^{'}$, $f^{''}$, $f^{'''}$ are first, second, and third derivatives of $f$ with respect to $t$. It  is invariant under all linear fractional transformations, namely $ \{ \frac{a f + b}{c f +d} , t \} = \{f , t\}$,
 where $ad - bc \neq 0$.  A geodesic $\gamma :I \rightarrow M$ is said to be projective if its natural parameter on $I$ is a projective parameter.
\section{Projectively invariant intrinsic distance  in complete Einstein spaces}
Consider the Funk metric $L_f$ defined on the open interval $I=\{u \in \mathbb{R} \mid -1 <u<1\}$ by
 $L_f = \frac{1}{k} (\frac{\mid y \mid}{1-u^2} + \frac{uy}{1-u^2})$, where $k$ is a constant. The Funk distance of every two points $a$ , $b$ $\in I$ is given by
\begin{equation}\label{e12}
D_f(a , b) = \frac{1}{2k} (\mid \ln \frac{(1-a)(1+b)}{(1-b)(1+a)} \mid + \ln \frac{(1-a^2)}{(1-b^2)}).
\end{equation}
See Refs.\cite{O} and \cite{SB} for more details.
Now, we are in a position to  define the pseudo-distance $d_M$, on a Finsler space $(M,F)$.
Given any two points $x,y \in M$, we choose a chain $\alpha$ of geodesic segments consisting of;
(I) a chain of points $x = x_0 , x_1 , ... ,x_k = y$ on $M$; (II) pairs of points $a_1,b_1 ,..., a_k,b_k$ in $I$; (III) projective maps $f_1,...,f_k$, $f_i: I \rightarrow M $, such that
$f_i(a_i) = x_{i-1}, \quad f_i(b_i) = x_i, \quad i = 1,...,k$. The length $L(\alpha)$ of the chain $\alpha$ is defined to be $L(\alpha) = \Sigma_i D_f(a_i , b_i)$. The pseudo-distance $d_M (x , y)$ is defined by $d_M(x , y) = inf L(\alpha)$, where the infimum is taken over all chains $\alpha$ from $x$ to $y$. It is well known $d_M$ remains  invariant under the projective change of metrics and we have the following lemmas, cf.,  \cite{SB}.
 \setcounter{thm}{0}
 \begin{lemma}\label{Lem;1}
(I) Let the geodesic  $f : I \rightarrow M $ be a projective map, then $D_f(a , b) \geq d_M (f(a) , f(b))$ for any $a,b \in I$.
(II) Let $\delta_M$ be any pseudo-distance on $M$ with the property  $D_f(a , b) \geq \delta_M(f(a) , f(b))$ for any $a,b \in I$, and for all projective maps $f: I \rightarrow M$, then
  $\delta_M(x , y)\leq d_M(x , y)$ for any $x,y \in M$.
  \end{lemma}
\begin{lemma}\label{Lem;2}
  Let $(M,F)$ be a Finsler space for which the Ricci tensor satisfies
 $ Ric_{ij} \leq -c^2g_{ij}$  as matrices, for a positive constant $c$. Let $d_F(. , .)$ be the distance induced by $F$, then for every projective map $f:I \rightarrow M$, $d_F$ is bounded below by the Funk distance $D_f$, that is, $D_f(a , b) \geq \frac{2 c}{\sqrt{n - 1} k} d_F(f(a) , f(b)) \quad \forall a,b \in I$.
 \end{lemma}
 \setcounter{thm}{0}
\begin{proposition}\label{Pro;1}
 Let $(M,F)$ be a Finsler space for which the Ricci tensor satisfies
 $ Ric_{ij} \leq -c^2g_{ij}$  as matrices, for a positive constant $c$. Then the pseudo-distance $d_M$ is a distance.
\end{proposition}
Following the procedure described above by collecting properties of the projectively invariant distance $d_M$,  we are in a position to  prove the following  theorem.
\setcounter{thm}{0}
\begin{theorem}\label{The;1}
 Let $(M , F)$ be a complete Einstein Finsler space with
\begin{equation}\label{e18}
(Ric)_{ij} = -c^2 g_{ij},
\end{equation}
 where $c$ is a positive constant. Then $d_M$ the projectively invariant distance is proportional to the Finslerian distance $d_F$, that is
 \begin{equation}\label{e19}
 d_M(x , y) = \frac{2c}{\sqrt{n - 1} k} d_F(x , y).
 \end{equation}
 \end{theorem}
 {\emph{Proof.}}
   By means of the second part of Lemma \ref{Lem;1}  and Lemma \ref{Lem;2}, we have
 $d_F(x , y) \frac{2c}{\sqrt{n - 1} k} \leq d_M(x , y)$.
  To prove the assertion, it remains to show  the inverse.
  Given any two points $x , y$ on $M$, we take a minimizing geodesic $x(s)$ on $M$ parameterized by  arc-length $s$ satisfying $x = x(0)$ and $y = x(a)$,  where $a$ is the Finslerian distance from $x$ to $y$.
  A projective parameter $\pi$ for this geodesic is defined to be a solution of the differential equation
  $\{\pi,s\}= \frac{2}{n-1} Ric_{jk}\frac{d{x}^j}{ds}\frac{d{x}^k}{ds}$.\\
  Let us consider the assumption (\ref{e18}) which leads to
  $\{\pi,s\}= \frac{- 2 c^2}{n-1} g_{jk}\frac{d{x}^j}{ds}\frac{d{x}^k}{ds}$.
For all curves parameterized by arc length $g_{jk}\frac{d{x}^j}{ds}\frac{d{x}^k}{ds} = 1$, therefore
  $\{\pi,s\}= \frac{- 2 c^2}{n-1}$.  General solution of the latter equation is Given by
  \begin{equation}\label{e20}
  \pi(s) = \frac{\alpha e^{js} + \beta e^{-js}}{\gamma e^{js} + \delta e^{-js}},
  \end{equation}
  where $\alpha \delta - \beta \gamma \neq 0$, and $j = \frac{c}{\sqrt{n - 1}}$.
 According to the first part of Lemma \ref{Lem;1}, the Funk distance $D_F$ between the points $0$ and $\pi{a}$ in $I$, satisfies
 \begin{equation}\label{e21}
  D_F(0 ,\pi(a)) \geq d_M(x , y).
  \end{equation}
   We consider the special solution of the equation (\ref{e20}) as follows $\pi(s) = \frac{e^{js} + e^{-js}}{e^{js} } = 1 - e^{-2js}$.
The real number $a$ is non-negative, hence $\pi(a)$ belongs to  the interval $I$. Thus plugging $\pi(s)$ into the equation (\ref{e12}) leads to
\begin{eqnarray*}
 D_f(0 , \pi(a))&=&\frac{1}{2k} \mid \ln (\frac{1 +\pi(a)}{1 - \pi(a)}) \mid+\frac{1}{2k} \ln \frac{1}{1 - (\pi(a))^2}=\frac{1}{2k} \mid \ln (\frac{2-e^{-2ja}}{e^{-2ja}}) \mid+ \frac{1}{2k} \ln \frac{e^{2ja}}{2 - e^{-2ja}} \\
  &=&\frac{1}{2k} \mid \ln (2 e^{2ja} - 1)\mid + \frac{1}{2k} \ln \frac{e^{2ja}}{2 - e^{-2ja}}.
\end{eqnarray*}
 By examining the fact $(2 e^{2ja} - 1) > 1$, we can omit the absolute value, $\mid . \mid$, which in turn leads to
  $$D_f(0 ,\pi(a))= \frac{1}{2k} \ln (2 e^{2ja} - 1)  + \frac{1}{2k}\ln e^{2ja}-\frac{1}{2k}\ln(2-e^{-2ja})$$
  $$= \frac{1}{2k} \ln (2 e^{2ja} - 1)+ \frac{ja}{k} - \frac{1}{2k}\ln \frac{2 e^{2ja} - 1}{e^{2ja}}= \frac{2ja}{k}=\frac{2cd_F(x,y)}{k\sqrt{n - 1}}.$$
   By plugging this relation into  the equation (\ref{e21}) we have  $d_M(x , y) \leq \frac{2c }{k \sqrt{n - 1}}d_F(x , y)$. This completes the proof.
  \hspace{\stretch{1}}$\Box$\\
  Two Finsler structures $F$ and $\bar{F}$ are said to be homothetic if there is a constant $\lambda$ such that  $F = \lambda \bar{F}$.
  A  result due to Busemann-Mayer theorem  for positively homogeneous functions shows that whenever we have a Finsler distance on a manifold $M$, we can find  the related Finsler structure, cf.,  \cite{BM} or \cite{BCS}, p.161.
   Considering the later assertion and Theorem \ref{The;1}, Corollary \ref{Cor;1} is easily obtained.
 \setcounter{thm}{0}
\begin{corollary}\label{Cor;1}
 Let $(M , F) $ and $(M , \bar{F})$ be two  complete Einstein Finsler spaces with $Ric_{ij} = -c^2 g_{ij}$ and  ${\bar{Ric}}_{ij} = -\bar{c}^2 {\bar{g}_{ij}}$ respectively, if  $F$ , $\bar{F}$ are projectively related then  they are homothetic.
\end{corollary}
Moreover, if the flag curvature $\lambda$ is constant then by means of the equation (\ref{ric,k}) we have the following corollary.
\begin{corollary}Let $(M , F) $ and $(M , \bar{F})$ be two complete Finsler spaces of constant negative flag curvature, if  $F$ and $\bar{F}$ are projectively related then they are homothetic.
\end{corollary}
 \ \ \ \ \

\end{document}